\documentclass[12pt]{amsart}
\begin{document}
\renewcommand{\thefootnote}{\fnsymbol{footnote}}
\begin{center}
{\bf\large A short constructive proof of Jordan's decomposition
theorem}
\medskip\\ {Pawel Kr\"oger\footnote{Research partially
supported by UTFSM Grant \# 12.02.23\\ Departamento de
Matem\'atica, Universidad T\'ecnica Federico Santa Mar\'{\i}a,
Valpara\'{\i}so, Chile\\e-mail: pawel.kroeger@mat.utfsm.cl}}
\end{center}
\medskip {\bf Abstract. }{\small Although there are many simple
proofs of Jordan's decomposition theorem in the literature (see
the bibliography and the references mentioned in $[$3$]$), our
proof seems to be even more elementary. In fact, all we need is
the theorem on the dimensions of range and kernel and the
existence of eigenvalues of a linear transformation on a
nontrivial finite dimensional complex vector space. We construct
whole chains and and count the total number of vectors belonging
to the chains. The proof that those vectors are linear independent
is standard. The Cayley-Hamilton theorem is not needed for the
proof, it could be obtained as an easy corollary.}
\bigskip\\
{\bf Theorem. }{\it Let $A$ be a linear operator on a finite
dimensional complex vector space $X$. There exist complex numbers
$\lambda_i$ and a base $\{{\bf x}_i^j\}$ of $X$ where $1\le i\le
n$ and $0\le j\le k_i$ such that $$A{\bf x}_i^0\,=\,\lambda_i{\bf
x}_i^0\quad\mbox{and}\quad A{\bf x}_i^j\,=\,\lambda_i{\bf
x}_i^j+{\bf x}_i^{j-1}\quad\mbox{for every }j\ge1\mbox{ and every
} i.\qquad(1)$$ \quad A second linear operator $B$ on $X$ is
similar to $A$ if and only if\\
$\dim(B-\lambda)^kX=\dim(A-\lambda)^kX$ for every eigenvalue
$\lambda$ and every $k\ge1$.}
\medskip\\ {\it Proof. } We prove the first part of the theorem by
induction with respect to the number of distinct eigenvalues of
$A$. If there are no eigenvalues, then the dimension of $X$ is $0$
and the theorem is trivial. Assume that the assertion is true for
every linear transformation with $d$ distinct eigenvalues.
Consider a linear transformation $A$ with $d+1$ distinct
eigenvalues. Let $\lambda$ be one of them. Set
$R_k=(A-\lambda)^kX$ for $k\ge0$. There exists a positive integer
$a$ such that $R_a$ is a proper subspace of $R_{a-1}$ and such
that $R_{a+1}=R_a$. Thus, the restriction of $A-\lambda$ to $R_a$
is invertible. The restriction of $A$ to $R_a$ has only $d$
distinct eigenvalues and the induction assumption applies. There
is a base $\{{\bf x}_i^j\}$ of $R_a$ where $1\le i\le n$ and $0\le
j\le k_i$ that satisfies (1). We need to extend that base to a
base of $X$.

The intersection of $R_a$ and the kernel of $(A-\lambda)^a$ is
$\{{\bf 0}\}$. Let $r_k$ be the dimension of $R_k$ for every $k$.
Let $N_k$ be the kernel of the restriction of $A-\lambda$ to $R_k$
for every $k$. Set $n_k=\dim N_k$. Thus, $n_k\ge n_{k+1}$ for
every $k$. By the theorem on the dimensions of range and kernel of
a linear transformation, $\dim N_k=r_k-r_{k+1}$.

For every $k$ with $n_k>n_{k+1}$ we will construct $n_k-n_{k+1}$
chains of length $k+1$ with the chain property (1) for the
eigenvalue $\lambda$. We choose inductively a base $\{{\bf
x}_i^0\}$ of $N_k$ for every $k\ge0$; here $n<i\le n+n_k$. For
$k\ge a$ we have $N_k=\{{\bf 0}\}, n_k=0$, and our base is empty.
Assume that we have chosen a base of $N_k$ for some $k$ with $1\le
k\le a$. We add $n_{k-1}-n_k$ vectors ${\bf x}_i^0$ for
$n+n_k<i\le n+n_{k-1}$ in order to obtain a base of $N_{k-1}$.

Now recall that $N_k\subset R_k$ for every $k$. Choose ${\bf
x}_i^k$ with $(A-\lambda)^k{\bf x}_i^k={\bf x}_i^0$ for
$n+n_{k+1}<i\le n+n_k$ and set ${\bf x}_i^j=(A-\lambda)^{k-j}{\bf
x}_i^k$ for $1\le j\le k-1$. Thus, $(A-\lambda)^{j+1}{\bf
x}_i^j=(A-\lambda){\bf x}_i^0={\bf 0}$ for every $i, j$ with
$i>n$. We obtain a total of $$\sum_{k=0}^{a-1}(k+1)(n_k-n_{k+1})
=\sum_{k=0}^{a-1}n_k=\sum_{k=0}^{a-1}(r_k-r_{k+1}) =\dim X-\dim
R_a$$ vectors from the kernel of $(A-\lambda)^a$. It remains to
show that those vectors are linear independent. We set $k_i=k$ for
$n+n_{k+1}<i\le n+n_k$ and $k<a$. Assume that there are complex
numbers $c_i^j$ such that $$\sum_{i=n+1}^{n+n_0}\sum_{j\le
k_i}c_i^j{\bf x}_i^j\,=\,{\bf 0}.$$ We prove by induction on $j$
that $c_i^j=0$ for every $i>n$ and every $j$. Recall that $k_i<a$
for every $i>n$. Assume that $c_i^j=0$ for every $j$ with $j>k$
and some $k\le a-1$. We apply $(A-\lambda)^k$ to the above
equation and obtain that $$\sum_{i=n+1}^{n+n_0}\sum_{j\le
k_i}c_i^j(A-\lambda)^k{\bf
x}_i^j\,=\,\sum_{i=n+1}^{n+n_k}c_i^k{\bf x}_i^0\,=\,{\bf 0}.$$
Thus, $c_i^k=0$ for $n<i\le n+n_k$ (recall that the vectors ${\bf
x}_i^0$ for $n<i\le n+n_k$ are linear independent).

The second part of the theorem is an obvious consequence of the
above construction.
\medskip\\
{\bf Acknowledgement. } The author is very grateful to Professor
Rainer Schimming and to Professor Mark Ashbaugh for bringing the
references $[$1$]$--$[$4$]$ to his attention.
\medskip\\
{\small {\bf References}
\medskip\\$[$1$]$ S. Cater. An elementary development of the Jordan
canonical form.\\ Amer. Math. Monthly 69 (1962), 391--393\\$[$2$]$
A. Galperin and Z. Waksman. An elementary approach to Jordan
theory.\\ Amer. Math. Monthly 87 (1981), 728--732\\ $[$3$]$ H.
V\"aliaho. An elementary approach to the Jordan form of a
matrix.\\ Amer. Math. Monthly 93 (1986), 711-714 \\ $[$4$]$ I.
Gohberg and S. Goldberg. A simple proof of the Jordan
decomposition theorem.\\ Amer. Math. Monthly 103 (1996), 157-159\\
$[$5$]$ P. Halmos. Finite-dimensional Vector Spaces. Undergraduate
Texts in Mathematics. Springer, 1996 (second edition)}
\end{document}